\documentclass[english]{amsart}
\usepackage{amssymb}

\providecommand{\tabularnewline}{\\}

\newcommand{\esssup}{\mathop{\mathrm{ess\,sup}}}

 \theoremstyle{plain}
\newtheorem{theorem}{Theorem}[section]
  \theoremstyle{plain}
  \newtheorem{corollary}[theorem]{Corollary}
  \theoremstyle{plain}
  \newtheorem{proposition}[theorem]{Proposition}
  \theoremstyle{plain}
  \newtheorem{lemma}[theorem]{Lemma}
  \theoremstyle{plain}
  \newtheorem{remark}[theorem]{Remark}

\begin{document}

\title[A Sampling Inequality]{A Sampling Inequality for Fractional Order Sobolev Semi-Norms Using Arbitrary Order Data}

\author{Andrew Corrigan}

\address{Andrew Corrigan \\ Department of Computational and Data Sciences \\ George Mason University\\ Fairfax, VA 22030, USA}

\author{John Wallin}

\address{John Wallin \\ Department of Computational and Data Sciences \\ George Mason University\\ Fairfax, VA 22030, USA}

\author{Thomas Wanner}

\address{Thomas Wanner \\ Department of Mathematical Sciences\\ George Mason University\\ Fairfax, VA 22030, USA}

\begin{abstract}
To improve convergence results obtained using a framework
for unsymmetric meshless methods due to Schaback~(Preprint G\"ottingen 2006), we
extend, in two directions, the Sobolev bound due to Arcang\'eli et al.~(Numer Math 107,
181-211, 2007), which itself extends two others due to Wendland and Rieger~(Numer Math 101,
643-662, 2005) and Madych~(J. Approx Theory 142, 116-128, 2006).
 The first is to incorporate discrete
samples of arbitrary order derivatives into the bound, which are used
to obtain higher order convergence in higher order Sobolev norms.
The second is to optimally bound fractional order Sobolev semi-norms,
which are used to obtain more optimal convergence rates when solving
problems requiring fractional order Sobolev spaces, notably inhomogeneous
boundary value problems.

\end{abstract}

\maketitle

\section{Introduction}

Over the past few years, increasingly general bounds of Sobolev
semi-norms, in terms of discrete samples, have appeared.  Such bounds are often 
called \emph{sampling inequalities}. A rather general sampling inequality was established
by Arcang\'{e}li et al.~\cite{arcangeli:2007:sampling_inequality}, and is stated as Theorem~\ref{thm:general_sampling_inequality}, with notation given in Section~\ref{sub:notation}.
\begin{theorem}
\label{thm:general_sampling_inequality}~\cite[Theorem 4.1]{arcangeli:2007:sampling_inequality}
Let $\Omega$ be a Lipschitz domain in $\mathbb{R}^{n}$, so that the domain $\Omega$ satisfies the cone property~\cite[Page 185]{arcangeli:2007:sampling_inequality}
 with radius $\rho>0$ and angle $\theta\in\left(0,\pi/2\right]$. Furthermore, let
  $p,q,\varkappa\in\left[1,\infty\right]$ and let $r$ be a real
number such that $r\geq n$, if $p=1$, $r>n/p$ if $1<p<\infty$,
or $r\in\mathbb{N}^{*}$, if $p=\infty$. Let $l_{0}=r-n\left(1/p-1/q\right)_{+}$
and $\gamma=\max\left\{ p,q,\varkappa\right\} $. Then, there exist
two positive constants $\mathfrak{d}_{r}$\textup{ (}\textup{\emph{dependent
on }}\textup{$\theta,\rho,n$ and $r$)} and $C$\textup{ (}\textup{\emph{dependent
on}}\textup{ $\Omega,n,r,p,q$ and $\varkappa$)} satisfying the following
property: for any set $A\subset\overline{\Omega}$ (or $A\subset\Omega$
if $p=1$ and $r=n$) such that $d=\delta\left(A,\overline{\Omega}\right)\leq\mathfrak{d}_{r}$ (c.f.~\eqref{eq:the_fill_distance})
for any $u\in W^{r,p}\left(\Omega\right)$ and for any real number
$l$ satisfying $l=0,\ldots,l_{max}$, we have \begin{equation}
\left|u\right|_{l,q,\Omega}\leq C\left(d^{r-l-n\left(1/p-1/q\right)_{+}}\left|u\right|_{r,p,\Omega}+d^{n/\gamma-l}\left\Vert \left.u\right|_{b}\right\Vert _{\varkappa}\right).\label{eq:arcangeli_sampling_inequality}\end{equation}
where $l_{\max}:=\left\lceil l_{0}\right\rceil -1$, unless the following
additional conditions hold, in which case $l_{\max}:=l_{0}$: $r\in\mathbb{N}^{*}$
and either (i) $p<q<\infty$ and $l_{0}\in\mathbb{N}$, (ii) $\left(p,q\right)=\left(1,\infty\right)$,
or (iii) $p\geq q$.
\end{theorem}
This sampling inequality generalizes those of Madych~\cite{madych:2006:sampling_inequality}
and Wendland and Rieger~\cite{wendland:2005:sampling_inequality},
by greatly extending the range of parameters $r,p,l,$ and $\varkappa$.
While Theorem~\ref{thm:general_sampling_inequality}
applies to functions with finite smoothness, an analogous bound for
functions with infinite smoothness has been provided by Rieger and
Zwicknagl~\cite{rieger:2006:exponential_sampling_inequality} which
achieves exponential factors.

Arcang\'{e}li, et al.~\cite{arcangeli:2007:sampling_inequality}
used Theorem~\ref{thm:general_sampling_inequality} to derive error
bounds for interpolating and smoothing $\left(m,s\right)$-splines,
an application which we do not consider. Instead we are interested
in another major application of these Sobolev estimates: Schaback's
framework for unsymmetric meshless methods for operator equations~\cite{schaback:2007:framework}, see also the earlier version~\cite{schaback:2007:meshless_collocation_convergence}.
A sampling inequality is
necessary for unsymmetric meshless methods, such as Schaback's modification
of Kansa's method~\cite{schaback:2007:framework,schaback:2007:meshless_collocation_convergence},
which involve an overdetermined system of equations in general. In
an attempt to improve the order of convergence obtained using Schaback's
framework, we extend the bound of Arcang\'{e}li, et al. in two ways.

Our first extension is to loosen the restriction $l\in\mathbb{N}$ to allow for fractional
order Sobolev norms on the left hand side of the sampling inequality.  In the context of Schaback's framework,
this will result in more optimal convergence results in terms of both the test and trial discretization parameters.   Otherwise,
the test discretization would require a higher rate of refinement.

Our second extension is to incorporate discrete samples of
arbitrary order derivatives into the bound. The reason for this is that~\eqref{eq:arcangeli_sampling_inequality} 
has a factor $d^{-l}$ in its second term which is insufficient for achieving a \emph{uniformly stable test discretization}
for higher order Sobolev norms in Schaback's framework. With this modification to incorporate
samples of higher order derivatives we will be able to come closer
to achieving such a test discretization, resulting in higher order
convergence results.  This introduces a new parameter $\mu$,
for which previous sampling inequalities coincide with the choice $\mu=0$.

\subsection{\label{sub:notation}Notation}

We employ the notation of Arcang\'{e}li, et al.~\cite[Section 2]{arcangeli:2007:sampling_inequality}.

As in~\cite[Section 4]{arcangeli:2007:sampling_inequality}, we assume throughout this paper
that $\Omega$ is a bounded domain in $\mathbb{R}^{n}$ with a Lipschitz-continuous
boundary, so that the domain $\Omega$ satisfies the cone property~\cite[Page 185]{arcangeli:2007:sampling_inequality} with radius $\rho>0$
and angle $\theta\in\left(0,\pi/2\right]$. For a given finite subset
$A$ of $\overline{\Omega}$, the \emph{fill distance} is defined
as 
\begin{equation}\label{eq:the_fill_distance}
\delta\left(A,\overline{\Omega}\right)=\sup_{x\in\Omega}\min_{a\in A}\left|x-a\right|.\end{equation}

The following restates a portion of their notation.
For all $r\in\left[0,\infty\right]$ and $p\in\left[1,\infty\right]$,
the Sobolev norm is denoted by  $\left\Vert \cdot\right\Vert _{r,p,\Omega}$,
while the Sobolev semi-norm is denoted by $\left|\cdot\right|_{r,p,\Omega}$.
The set $\mathbb{N}^{*}=\left\{1,2,3,\ldots\right\}$, while $\mathbb{N}=\left\{0,1,2,\ldots\right\}$.
The space of polynomials over 
$\mathbb{R}^n$ with degree less than or equal to $k$ is denoted by $P_{k}$.

We make the following additions to their notation.
Let $H^r\left(\Omega\right):= W^{r,2}\left(\Omega\right)$
and \[
\tilde{W}^{r,q}\left(\Omega\right):=\left\{ v\in W^{r,q}\left(\Omega\right)\,:\,\int_{\Omega}v=0\right\} .\]
Given a function $v\in W^{1,q}\left(\Omega\right)$, the vector-valued
function consisting of its partial derivatives is denoted by $Dv$.
The surface area of the $n$-dimensional
ball is denoted by $\left|S^{n-1}\right|$. In Section~\ref{sec:Preliminaries}, a generic constant $C$ appears in many proofs, whose particular
value may change, but with the parameters on which it depends either
indicated in parentheses or stated explicitly in the exposition. We
will often substitute dependencies with others, possibly taking the
maximum or minimum value, as required by its application, of the constant
over a finite range of values. 
In Section~\ref{sec:Unsymmetric-Meshless-Methods}, only the dependence of constants on the discretization
parameters $r$ and $s$ is explicitly stated since we regard the
spaces and mappings in that section as fixed.

\section{\label{sec:Preliminaries}Extension of the Sobolev Bound}

\subsection{\label{sub:Sobolev-spaces}Fractional Order Sobolev Spaces}

This section concerns fractional order Sobolev norms and the results
of this section will be used to generalize~\cite[Proposition 3.4]{arcangeli:2007:sampling_inequality}
to Proposition~\ref{pro:ball_high_order_sampling_inequality}. Lemmas
\ref{lem:mean_value_Lq} and~\ref{lem:mean_value_Linf} each require
an extension operator which satisfies~\eqref{eq:seminorm_extension_bound}
for zero-average functions over a ball. This property is not provided
by standard extension operators since they involve a domain-dependent
constant and the full Sobolev norm in the bound, rather than a domain-independent
constant and the Sobolev semi-norm.
\begin{lemma}
\label{lem:seminorm_extension}If $q\in\left[1,\infty\right]$, $r>0$,
and $x_{0}\in\mathbb{R}^{n}$ then there exists a linear, continuous
operator \[
E:\tilde{W}^{1,q}\left(B\left(x_{0},r\right)\right)\rightarrow W^{1,q}\left(\mathbb{R}^{n}\right)\]
 such that for all $v\in\tilde{W}^{1,q}\left(B\left(x_{0},r\right)\right)$\begin{equation}
Ev=v\mbox{ a.e. in }B\left(x_{0},r\right),\label{eq:seminorm_extension_property}\end{equation}
and \begin{equation}
\left|Ev\right|_{1,q,\mathbb{R}^{n}}\leq C\left(n,q\right)\left|v\right|_{1,q,B\left(x_{0},r\right)}.\label{eq:seminorm_extension_bound}\end{equation}
If in addition $v\in C^{1}\left(\overline{B\left(x_{0},r\right)}\right)$
then $Ev\in C^{1}\left(\mathbb{R}^{n}\right)$.\end{lemma}
\begin{proof}
Let $v\in\tilde{W}^{1,q}\left(B\left(x_{0},r\right)\right)$ and $\hat{v}:=v\circ F$
where $F:\hat{x}\rightarrow r\hat{x}+x_{0}$. From a change of variables
it follows that $\hat{v}\in W^{1,q}\left(B\left(0,1\right)\right)$
with semi-norm \[
\left|\hat{v}\right|_{1,q,B\left(0,1\right)}=r^{1-n/q}\left|v\right|_{1,q,B\left(x_{0},r\right)}.\]
From~\cite[Section 5.4, Theorem 1]{evans:pde}, there exists a linear,
continuous extension operator \[
\hat{E}:W^{1,q}\left(B\left(0,1\right)\right)\rightarrow W^{1,q}\left(\mathbb{R}^{n}\right)\]
such that for each $v\in W^{1,q}\left(B\left(0,1\right)\right)$\[
\hat{E}\hat{v}=\hat{v}\mbox{ a.e. in }B\left(0,1\right),\]
and\[
\left\Vert \hat{E}\hat{v}\right\Vert _{1,q,\mathbb{R}^{n}}\leq C\left(n,q\right)\left\Vert \hat{v}\right\Vert _{1,q,B\left(0,1\right)},\]
where the dependence on $n$ is through $B\left(0,1\right)$. That \begin{eqnarray*}
\left\Vert \hat{v}\right\Vert _{1,q,B\left(0,1\right)} & \leq & C\left(n,q\right)\left|\hat{v}\right|_{1,q,B\left(0,1\right)},\end{eqnarray*}
follows from specializing a Poincar\'{e} inequality given in~\cite[Section 5.8, Theorem 2]{evans:pde}
to the unit ball and that \begin{equation}\int_{B\left(0,1\right)}\hat{v}=r^{-n}\int_{B\left(x_{0},r\right)}v=0.\end{equation}
Let \begin{equation}
Ev:=\left(\hat{E}\hat{v}\right)\circ F^{-1}\label{eq:definition_of_extension}\end{equation}
so that the result~\eqref{eq:seminorm_extension_property} holds.
From another change of variables, it follows that \[
\left|Ev\right|_{1,q,\mathbb{R}^{n}}=r^{n/q-1}\left|\hat{E}\hat{v}\right|_{1,q,\mathbb{R}^{n}}.\]
The result~\eqref{eq:seminorm_extension_bound} follows by combining
the preceding relations. Finally, from the proof of~\cite[Section 5.4, Theorem 1]{evans:pde}
it follows that if $\hat{v}\in C^{1}\left(\overline{B\left(0,1\right)}\right),$
then $\hat{E}\hat{v}\in C^{1}\left(\mathbb{R}^{n}\right)$, so that
if $v\in C^{1}\left(\overline{B\left(x_{0},r\right)}\right)$ then
$Ev\in C^{1}\left(\mathbb{R}^{n}\right)$.
\end{proof}
Based on this extension, we obtain Lemma~\ref{lem:mean_value_Lq},
which is similar to a result used by Bourgain et al.~\cite[Eq. 2]{book_with_fractional_paper:2001},
but uses a domain-independent constant and semi-norm in the bound.
\begin{lemma}
\label{lem:mean_value_Lq}If $q\in\left[1,\infty\right)$, $h\in\mathbb{R}^{n}$,
and $v\in\tilde{W}^{1,q}\left(B\left(x_{0},r\right)\right)$ then\[
\left(\int_{B\left(x_{0},r\right)}\left|Ev\left(x+h\right)-Ev\left(x\right)\right|^{q}dx\right)^{1/q}\leq C\left(n,q\right)\left|h\right|\left|v\right|_{1,q,B\left(x_{0},r\right)}\]
\end{lemma}
\begin{proof}
First suppose that $v$ is in the subset \begin{equation}C^{1}\left(\overline{B\left(x_{0},r\right)}\right)\cap\tilde{W}^{1,q}\left(B\left(x_{0},r\right)\right)\end{equation}
which is dense in $\tilde{W}^{1,q}\left(B\left(x_{0},r\right)\right)$.
From Lemma~\ref{lem:seminorm_extension} it follows that \begin{eqnarray*}
\int_{B\left(x_{0},r\right)}\left|Ev\left(x+h\right)-Ev\left(x\right)\right|^{q}dx & \leq & \int_{\mathbb{R}^{n}}\left|Ev\left(x+h\right)-Ev\left(x\right)\right|^{q}dx\\
 & = & \int_{\mathbb{R}^{n}}\left|\int_{0}^{1}\frac{d}{dt}Ev\left(x+th\right)dt\right|^{q}dx\\
 & \leq & \int_{\mathbb{R}^{n}}\int_{0}^{1}\left|\frac{d}{dt}Ev\left(x+th\right)\right|^{q}dtdx\\
 & = & \int_{\mathbb{R}^{n}}\int_{0}^{1}\left|DEv\left(x+th\right)\cdot h\right|^{q}dtdx\\
 & \leq & \left|h\right|^{q}\int_{\mathbb{R}^{n}}\int_{0}^{1}\left|DEv\left(x+th\right)\right|^{q}dtdx\\
 & = & \left|h\right|^{q}\int_{0}^{1}\int_{\mathbb{R}^{n}}\left|DEv\left(x+th\right)\right|^{q}dxdt\\
 & \leq & \left|h\right|^{q}\int_{0}^{1}\left|Ev\left(\cdot+th\right)\right|_{1,q,\mathbb{R}^{n}}^{q}dt\\
 & = & \left|h\right|^{q}\left|Ev\right|_{1,q,\mathbb{R}^{n}}^{q}\\
 & \leq & C\left(n,q\right)\left|h\right|^{q}\left|v\right|_{1,q,B\left(x_{0},r\right)}^{q}\end{eqnarray*}
From this, the result follows for all functions in $\tilde{W}^{1,q}\left(B\left(x_{0},r\right)\right)$
via a standard density argument.\end{proof}
\begin{lemma}
\label{lem:mean_value_Linf}If $x,x+h\in B\left(x_{0},r\right)$,
and $v\in\tilde{W}^{1,\infty}\left(B\left(x_{0},r\right)\right)$
then\[
\left|v\left(x+h\right)-v\left(x\right)\right|\leq C\left(n\right)\left|h\right|\left|v\right|_{1,\infty,B\left(x_{0},r\right)}\]
\end{lemma}
\begin{proof}
In the proof of~\cite[Section 5.8.2, Theorem 4]{evans:pde} it is
shown that $v$ is a Lipschitz function with constant $\left|Ev\right|_{1,\infty,\mathbb{R}^{n}}$,
where the extension operator constructed in~\cite[Section 5.4, Theorem 1]{evans:pde}
is used. However, the extension operator from Lemma~\ref{lem:seminorm_extension}
could be substituted so that the result then follows by~\eqref{eq:seminorm_extension_bound}.
\end{proof}
In the bound provided by Lemma~\ref{lem:fractional_1_norm_relation}
the factor $r^{1-\epsilon}$ will be the key to generalizing sampling
inequalities to optimally bound fractional order semi-norms.
\begin{proposition}
\label{lem:fractional_1_norm_relation}If $q\in\left[1,\infty\right]$,
$\epsilon\in\left(0,1\right)$, and $v\in W^{1,q}\left(B\left(x_{0},r\right)\right)$
then \begin{equation}
\left|v\right|_{\epsilon,q,B\left(x_{0},r\right)}\leq C\left(n,q\right)\left(1-\epsilon\right)^{-1/q}r^{1-\epsilon}\left|v\right|_{1,q,B\left(x_{0},r\right)}.\label{eq:fractional_one_seminorm_relation_on_the_ball}\end{equation}
\end{proposition}
\begin{proof}
Since the semi-norms that appear in~\eqref{eq:fractional_one_seminorm_relation_on_the_ball}
are invariant with respect to a shift in value of $v$ by a constant,
it suffices to only consider $v\in\tilde{W}^{1,q}\left(B\left(x_{0},r\right)\right)$.

Case $q\in\left[1,\infty\right)$: Let $y\in B\left(x_{0},r\right)$
and \[
B\left(x_{0},r\right)-y:=\left\{ x-y\,:\, x\in B\left(x_{0},r\right)\right\} ,\]
so that $B\left(x_{0},r\right)-y\subseteq B\left(0,2r\right).$ 

\begin{eqnarray*}
\left|v\right|_{\epsilon,q,B\left(x_{0},r\right)}^{q} & = & \int_{B\left(x_{0},r\right)}\int_{B\left(x_{0},r\right)}\frac{\left|v\left(x\right)-v\left(y\right)\right|^{q}}{\left|x-y\right|^{n+\epsilon q}}dxdy\\
 & = & \int_{B\left(x_{0},r\right)}\int_{B\left(x_{0},r\right)-y}\frac{\left|v\left(y+h\right)-v\left(y\right)\right|^{q}}{\left|h\right|^{n+\epsilon q}}dhdy\\
 & \leq & \int_{B\left(x_{0},r\right)}\int_{B\left(0,2r\right)}\frac{\left|Ev\left(y+h\right)-Ev\left(y\right)\right|^{q}}{\left|h\right|^{n+\epsilon q}}dhdy\\
 & = & \int_{B\left(0,2r\right)}\frac{\int_{B\left(x_{0},r\right)}\left|Ev\left(y+h\right)-Ev\left(y\right)\right|^{q}dy}{\left|h\right|^{n+\epsilon q}}dh\\
 & \leq & C\left(n,q\right)\left(\int_{B\left(0,2r\right)}\frac{\left|h\right|^{q}}{\left|h\right|^{n+\epsilon q}}dh\right)\left|v\right|_{1,q,B\left(x_{0},r\right)}^{q}\\
 & = & C\left(n,q\right)\left|S^{n-1}\right|\left(\int_{0}^{2r}\frac{\rho^{q}}{\rho^{n+\epsilon q}}\rho^{n-1}d\rho\right)\left|v\right|_{1,q,B\left(x_{0},r\right)}^{q}\\
 & \leq & \frac{C\left(n,q\right)}{\left(1-\epsilon\right)q}\left|S^{n-1}\right|\left(2r\right)^{\left(1-\epsilon\right)q}\left|v\right|_{1,q,B\left(x_{0},r\right)}^{q}\\
 & \leq & \frac{C\left(n,q\right)}{\left(1-\epsilon\right)q}\left|S^{n-1}\right|2^{q}r^{\left(1-\epsilon\right)q}\left|v\right|_{1,q,B\left(x_{0},r\right)}^{q}.\end{eqnarray*}

Case $q=\infty$: \begin{eqnarray*}
\left|v\right|_{\epsilon,q,B\left(x_{0},r\right)} & = & \esssup_{x,y\in B\left(x_{0},r\right),x\neq y}\frac{\left|v\left(x\right)-v\left(y\right)\right|}{\left|x-y\right|^{\epsilon}}\\
 & \leq & C\left(n,q\right)\left(2r\right)^{1-\epsilon}\left|v\right|_{1,q,B\left(x_{0},r\right)}\\
 & \leq & 2C\left(n,q\right)r^{1-\epsilon}\left|v\right|_{1,q,B\left(x_{0},r\right)}.\end{eqnarray*}
\end{proof}
\begin{remark}
\label{rem:defect} The explicit constant $\left(1-\epsilon\right)^{-1/q}$,
which blows up as $\epsilon$ increases towards one, is a manifestation
of the ``defect'' of intrinsic fractional order Sobolev semi-norms
studied by Bourgain et al.~\cite{book_with_fractional_paper:2001}.\end{remark}
\begin{corollary}
\label{lem:fractional_ceil_norm_relation}If $q\in\left[1,\infty\right]$,
$l\in\left[0,\infty\right]$, and $v\in W^{\left\lceil l\right\rceil ,q}\left(B\left(x_{0},r\right)\right)$
then\begin{equation}
\left|v\right|_{l,q,B\left(x_{0},r\right)}\leq C\left(n,q,\left\lfloor l\right\rfloor \right)K\left(\left\lceil l\right\rceil -l,q\right)r^{\left\lceil l\right\rceil -l}\left|v\right|_{\left\lceil l\right\rceil ,q,B\left(x_{0},r\right)},\label{eq:bound_fractional_bound_by_ceil}\end{equation}
where \begin{equation}
K\left(\left\lceil l\right\rceil -l,q\right):=\left\{ \begin{array}{cc}
1 & \mathrm{for}\:l\in\mathbb{N}\mbox{ or }q=\infty\\
\left(\left\lceil l\right\rceil -l\right)^{-1/q} & \mathrm{for}\:l\notin\mathbb{N}\mbox{ and }q<\infty\end{array}\right..\label{eq:correction_factor}\end{equation}

\end{corollary}

\subsection{An Auxiliary Result}

The following result applies Corollary~\ref{lem:fractional_ceil_norm_relation}
to generalize~\cite[Proposition 3.4]{arcangeli:2007:sampling_inequality}.
\begin{proposition}
\label{pro:ball_high_order_sampling_inequality}Let $p,q,\varkappa\in\left[1,\infty\right]$
such that $p\leq q$. Let $r$ be a real number such that $r>n/p$,
if $p>1$, or $r\geq n$, if $p=1$. Finally, let $k=\left\lceil r\right\rceil -1$,
$\mathfrak{K}=\dim P_{k}$, and $l_{0}=r-n/p+n/q$. Then, there exists a
constant $R>1$ (dependent on $n$ and $r$) and, for any $M'\geq1$,
there exists two constants $C$ (dependent on \textup{$M',n,r,p,q,$
and $\varkappa$)} \textup{and $K\geq1$ (explicitly dependent on
$\left\lceil l\right\rceil -l$ and $q$, cf.~\eqref{eq:correction_factor}),}
satisfying the following property: for any $d>0$ and any $t\in\mathbb{R}^{n}$,
the open ball $B\left(t,Rd\right)$ contains $\mathfrak{K}$ closed balls
$\mathcal{B}_{1},\ldots\mathcal{B}_{\mathfrak{K}}$ of radius $d$ such
that, for any $v\in W^{r,p}\left(\overline{B}\left(t,M'Rd\right)\right)$,
for any $b\in\Pi_{i=1}^{\mathfrak{K}}\mathcal{B}_{i}$ and $l\in\left[0,l_{\max}\right]$,\begin{equation}
\left|v\right|_{l,q,\overline{B}\left(t,M'Rd\right)}\leq C\cdot K\left(d^{r-l-n/p+n/q}\left|v\right|_{r,p,\overline{B}\left(t,M'Rd\right)}+d^{n/q-l}\left\Vert \left.v\right|_{b}\right\Vert _{\varkappa}\right),\label{eq:sampling_inequality_over_ball}\end{equation}
where we have let $l_{\max}:=\left\lceil l_{0}\right\rceil -1$, or
$l_{\max}:=l_{0}$ if the following additional conditions hold: $r\in\mathbb{N}^{*}$
and either (i) $p<q<\infty$ and $l_{0}\in\mathbb{N}$, (ii) $\left(p,q\right)=\left(1,\infty\right)$,
or (iii) $1\leq p=q\leq\infty$. \end{proposition}
\begin{proof}
The case that $l\in\mathbb{N}$ is established by~\cite[Proposition 3.4]{arcangeli:2007:sampling_inequality}.
Suppose that $l\notin\mathbb{N}$. The hypotheses imply that $l_{\max}\in\mathbb{N}$,
so that $\left\lceil l\right\rceil \leq l_{\max}$, and thus the result
follows by combining~\eqref{eq:sampling_inequality_over_ball} for
$l=\left\lceil l\right\rceil $ with Corollary~\ref{lem:fractional_ceil_norm_relation},
using the fact that $\left(M'R\right)^{\left\lceil l\right\rceil -l}\leq M'R$,
and substituting the dependence on $\left\lfloor l\right\rfloor $
and $R$ with $n,r,p$ and $q$.
\end{proof}

\subsection{Sobolev Bounds}

For $p\leq q$, the following result generalizes~\cite[Theorem 4.1]{arcangeli:2007:sampling_inequality}
to bound fractional order Sobolev semi-norms. No generalization to
bound fractional order Sobolev semi-norms is made for $p>q$, since
we have not obtained the relation~\cite[Eq. 2.1]{arcangeli:2007:sampling_inequality}
for the case that $l$ is fractional.
\begin{theorem}
\label{thm:fractional_high_order_sampling_inequality}Let $p,q,\varkappa\in\left[1,\infty\right]$
and let $r$ be a real number and $\mu$ a nonnegative integer such
that $r-\mu\geq n$, if $p=1$, $r-\mu>n/p$ if $1<p<\infty$, or
$r-\mu\in\mathbb{N}^{*}$ if $p=\infty$. Let $l_{0}=r-\mu-n\left(1/p-1/q\right)_{+}$
and $\gamma=\max\left\{ p,q,\varkappa\right\} $. Then, there exist
three positive constants $\mathfrak{d}_{r}$\textup{ (dependent on $\theta,\rho,n,r$
and $\mu$),} $C$\textup{ (dependent on $\Omega,n,r,p,q,$ and $\varkappa$),
and $K\geq1$ (explicitly dependent on $\left\lceil l\right\rceil -l$
and $q$, cf.~\eqref{eq:correction_factor}),} satisfying the following
property: for any set $A\subset\overline{\Omega}$ (or $A\subset\Omega$
if $p=1$ and $r-\mu=n$) such that $d=\delta\left(A,\overline{\Omega}\right)\leq\mathfrak{d}_{r}$,
for any $u\in W^{r,p}\left(\Omega\right)$ and, if $p\leq q$ then
for any real number $l\in\left[0,l_{\max}\right]$, otherwise if $p>q$
then for any integer $l=0,\ldots,l_{\max}$, \begin{equation}
\left|u\right|_{l,q,\Omega}\leq C\cdot K\left(d^{r-l-n\left(1/p-1/q\right)_{+}}\left|u\right|_{r,p,\Omega}+d^{n/\gamma+\mu-l}\left\Vert \prod_{\left|\alpha\right|=\mu}\left.\partial^{\alpha}u\right|_{A}\right\Vert _{\varkappa}\right),\label{eq:penultimate_seminorm_sampling_inequality}\end{equation}
where we have let $l_{\max}:=\left\lceil l_{0}\right\rceil -1$, or
$l_{\max}:=l_{0}$ if the following additional conditions hold: $r\in\mathbb{N}^{*}$
and either (i) $p<q<\infty$ and $l_{0}\in\mathbb{N}$, (ii) $\left(p,q\right)=\left(1,\infty\right)$,
or (iii) $p\geq q$. \end{theorem}
\begin{proof}
The case that $\mu=0$ and $l\in\mathbb{N}$ is~\cite[Theorem 4.1]{arcangeli:2007:sampling_inequality}.
The proof of this theorem for $\mu=0$ and $l\notin\mathbb{N}$ can
be obtained by reusing the proof of~\cite[Theorem 4.1]{arcangeli:2007:sampling_inequality},
but applying Proposition~\ref{pro:ball_high_order_sampling_inequality}
instead of~\cite[Proposition 3.4]{arcangeli:2007:sampling_inequality},
which allows for $l$ to be of fractional order for $p\leq q$, and
introduces the constant $K$. We now consider the case $\mu>0$. Let
$\alpha$ be a multi-index such that $\left|\alpha\right|=\mu$ and
therefore $\partial^{\alpha}u\in W^{r-\mu,p}\left(\Omega\right)$.
It follows from the case that $\mu=0$ that in the situation required
by the present hypotheses \[
\left|\partial^{\alpha}u\right|_{l-\mu,q,\Omega}\leq C\cdot K\left(d^{r-l-n\left(1/p-1/q\right)_{+}}\left|\partial^{\alpha}u\right|_{r-\mu,p,\Omega}+d^{n/\gamma+\mu-l}\left\Vert \left.\partial^{\alpha}u\right|_{A}\right\Vert _{\varkappa}\right).\]
We have for all $\alpha$ satisfying $\left|\alpha\right|=\mu$ that\[
\left|\partial^{\alpha}u\right|_{r-\mu,p,\Omega}\leq\left|u\right|_{r,p,\Omega}\]
and\[
\left\Vert \left.\partial^{\alpha}u\right|_{A}\right\Vert _{\varkappa}\leq\left\Vert \prod_{\left|\beta\right|=\mu}\left.\partial^{\beta}u\right|_{A}\right\Vert _{\varkappa}.\]
The result follows immediately for $q=\infty$. Otherwise, if $1\leq q<\infty$,
using that\[
\left|u\right|_{l,q,\Omega}^{q}\leq\sum_{\left|\alpha\right|=\mu}\left|\partial^{\alpha}u\right|_{l-\mu,q,\Omega}^{q}\]
the results follows with an additional factor $\left(\#\left\{ \alpha:\left|\alpha\right|=\mu\right\} \right)^{1/q}$
in the constant $C$, whose dependence on $\mu$ can be substituted
with $n,r$ and $p$.\end{proof}
\begin{corollary}
\label{cor:norm_sampling_inequality}Given the situation of Theorem
\ref{thm:fractional_high_order_sampling_inequality} with a constant
$\mathfrak{d}_{r}$ now dependent on \textup{$\theta,\rho,n,r,p$ and
$q$}, and the additional assumption that $r-l\in\mathbb{N}$, then we have \[
\left\Vert u\right\Vert _{l,q,\Omega}\leq C\cdot K\left(d^{r-l-n/p+n/q}\left\Vert u\right\Vert _{r,p,\Omega}+d^{n/\gamma+\mu-l}\left\Vert \prod_{\left|\alpha\right|\leq\mu}\left.\partial^{\alpha}u\right|_{A}\right\Vert _{\varkappa}\right).\]
\end{corollary}
\begin{proof}
From Theorem~\ref{thm:fractional_high_order_sampling_inequality},
there exists three positive constants $\mathfrak{d}_{r}\left(\theta,\rho,n,r,p,q\right)$,
$C\left(\Omega,n,r,p,q,\varkappa\right)$, and $K\left(\left\lceil l\right\rceil -l,q\right)\geq1$,
cf.~\eqref{eq:correction_factor},  such that for $d\leq\mathfrak{d}_{r}$
and $\eta=0,\ldots,\left\lfloor l\right\rfloor $ 
\begin{eqnarray*}
\left|u\right|_{\eta,q,\Omega}&\leq&C\cdot K\left(d^{r-l-n/p+n/q}\left|u\right|_{r-l+\eta,p,\Omega}\right.\\
& & ~~~~~~\left.+d^{n/\gamma+\left(\eta+\mu-\left\lfloor l\right\rfloor \right)_{+}-\eta}\left\Vert \prod_{\left|\alpha\right|=\left(\eta+\mu-\left\lfloor l\right\rfloor \right)_{+}}\left.\partial^{\alpha}u\right|_{A}\right\Vert _{\varkappa}\right)\label{eq:semi_sampling-inequality-for-each-eta}
\end{eqnarray*}
and for $\eta=l$ that
\begin{eqnarray*}
\left|u\right|_{l,q,\Omega} & \leq & C\cdot K\left(d^{r-l-n/p+n/q}\left|u\right|_{r,p,\Omega}+d^{n/\gamma+\mu-l}\left\Vert \prod_{\left|\alpha\right|=\mu}\left.\partial^{\alpha}u\right|_{A}\right\Vert _{\varkappa}\right).\label{eq:highest_order_semi_sampling_inequality}
\end{eqnarray*}
We have taken the constants to be the minimum or maximum over $\eta=0,\ldots,\left\lfloor l\right\rfloor ,l$
as required. Additionally we have restricted $\mathfrak{d}_{r}$ to be
at most one. For $\eta=0,\ldots,\left\lfloor l\right\rfloor $, we
have applied Theorem~\ref{thm:fractional_high_order_sampling_inequality}
with $r=r-l+\eta$ and $\mu=\left(\eta+\mu-\left\lfloor l\right\rfloor \right)_{+}$,
introducing a dependence of $\mathfrak{d}_{r}$ on $p$ and $q$ through
$l$, which along with $n$ and $r$ has substituted for the dependence
on $\mu$. It also follows for all $\eta=0,\ldots,\left\lfloor l\right\rfloor $
that\begin{equation}
\left\Vert \prod_{\left|\alpha\right|=\left(\eta+\mu-\left\lfloor l\right\rfloor \right)_{+}}\left.\partial^{\alpha}u\right|_{A}\right\Vert _{\varkappa}\leq\left\Vert \prod_{\left|\alpha\right|\leq\mu}\left.\partial^{\alpha}u\right|_{A}\right\Vert _{\varkappa},\label{eq:semi-samples-bounded-by-full-samples}\end{equation}
with a similar bound holding for $\left\Vert \prod_{\left|\alpha\right|=\mu}\left.\partial^{\alpha}u\right|_{A}\right\Vert _{\varkappa}$.
It follows from \[
\left(\eta+\mu-\left\lfloor l\right\rfloor \right)_{+}-\eta\geq\mu-\left\lfloor l\right\rfloor \geq\mu-l\]
and $\mathfrak{d}_{r}\leq1$ that for all $d\leq\mathfrak{d}_{r}$, \begin{equation}
d^{n/\gamma+\left(\eta+\mu-\left\lfloor l\right\rfloor \right)_{+}-\eta}\leq d^{n/\gamma+\mu-l}.\label{eq:samples-factor-bound}\end{equation}
It follows from $r-l\in\mathbb{N}$ that for each $\eta=0,\ldots,\left\lfloor l\right\rfloor $
we have $r-l+\eta\in\mathbb{N}$ and $0\leq r-l+\eta\leq r-l$, which
implies that \begin{equation}
\left|u\right|_{r-l+\eta,p,\Omega}\leq\left\Vert u\right\Vert _{r,p,\Omega}.\label{eq:semi_norm_by_full_norm}\end{equation}
Combining the preceding bounds we obtain for all $d\leq\mathfrak{d}_{r}$
and $\eta=0,\ldots,\left\lfloor l\right\rfloor ,l$ that \begin{equation}
\left|u\right|_{\eta,q,\Omega}\leq C\cdot K\left(d^{r-l-n/p+n/q}\left\Vert u\right\Vert _{r,p,\Omega}+d^{n/\gamma+\mu-l}\left\Vert \prod_{\left|\alpha\right|\leq\mu}\left.\partial^{\alpha}u\right|_{A}\right\Vert _{\varkappa}\right).\label{eq:semi-sampling-inequality-for-each-eta-up-to-fractional-order}\end{equation}
If $q=\infty$ then the result follows immediately. If $1\leq q<\infty$
it follows from~\eqref{eq:semi-sampling-inequality-for-each-eta-up-to-fractional-order}
that 
\begin{eqnarray*}
\left\Vert u\right\Vert _{l,q,\Omega} & \leq & C\cdot K\cdot\left(\left\lfloor l\right\rfloor +2\right)^{1/q}\left(d^{r-l-n/p+n/q}\left\Vert u\right\Vert_{r,p,\Omega}\right.\\
 &  & ~~~~~~~~~~~~~~~~~~~~~~~\left.+d^{n/\gamma+\mu-l}\left\Vert \prod_{\left|\alpha\right|\leq\mu}\left.\partial^{\alpha}u\right|_{A}\right\Vert _{\varkappa}\right).\end{eqnarray*}
The result then follows by incorporating the constant $\left(\left\lfloor l\right\rfloor +2\right)^{1/q}$
into $C$, with the dependence on $\left\lfloor l\right\rfloor $
substituted with $n,r,p$ and $q$.
\end{proof}
Only the case that $p=q=\varkappa=2$ will be used in Section~\ref{sec:Unsymmetric-Meshless-Methods}.
\begin{corollary}
\label{cor:specialized_sampling_inequality}Let $r$ be a real number,
$\mu$ be a nonnegative integer such that $r-\mu>n/2$. Then, there
exist three positive constants $\mathfrak{d}_{r}$\textup{ (dependent
on $\theta,\rho,n$ and $r$),} $C$\textup{ (dependent on $\Omega,n$
and $r$), and $K$ (explicitly dependent on $\left\lceil l\right\rceil -l$, cf.~\eqref{eq:correction_factor} with $q=2$)}
satisfying the following property: for any set $A\subset\overline{\Omega}$,
such that $d=\delta\left(A,\overline{\Omega}\right)\leq\mathfrak{d}_{r}$,
$u\in W^{r,2}\left(\Omega\right)$ and real number $l\in\left[0,r-\mu\right]$
such that $r-l\in\mathbb{N}$, \begin{equation}
\left\Vert u\right\Vert _{l,2,\Omega}\leq C\cdot K\left(d^{r-l}\left\Vert u\right\Vert _{r,2,\Omega}+d^{n/2+\mu-l}\left\Vert \prod_{\left|\alpha\right|\leq\mu}\left.\partial^{\alpha}u\right|_{A}\right\Vert _{2}\right).\label{eq:specialized_sampling_inequality}\end{equation}
 
\end{corollary}

\section{\label{sec:Unsymmetric-Meshless-Methods}Application: Unsymmetric
Meshless Methods for Operator Equations}

In this section, only the dependence of constants on the discretization
parameters $r$ and $s$ is explicitly stated since we regard the
spaces and mappings involved as fixed.

We now apply the sampling inequality stated in Corollary~\ref{cor:specialized_sampling_inequality} to
Schaback's framework for unsymmetric meshless methods for operator equations~\cite{schaback:2007:framework}.
Due to unaddressed issues contained
in its original formulation, we use a modified version stated in this section.
On a technical level, it differs substantially, e.g., certain spaces
have been eliminated, the inequalities apply over possibly different spaces, and the proof 
of the error bound has been slightly modified, but the underlying ideas are the same and 
entirely due to Schaback~\cite{schaback:2007:framework}.

The framework provides an error bound for meshless methods
which approximately solve a linear operator equation in the following setting.
The first requirement is a continuous and bijective linear operator
$L:U\rightarrow F$ mapping from the \emph{solution space} to the
\emph{data space}. The spaces $U$ and $F$ are assumed to be complete
in order to ensure the boundedness of $L^{-1}:F\rightarrow U$.
It is also assumed that the
exact solution $u^{*}\in \tilde{U}$ where $\tilde{U}\subset U$ is called the \emph{regularity subspace}. We denote
$\tilde{F}:=L\tilde{U}$.
The framework requires a scale of finite-dimensional trial subspaces $U_{r}\subset\tilde{U}$
equipped with a projector $\Pi_{r}:\tilde{U}\rightarrow U_{r}$. The
framework requires a linear, continuous, and bijective \emph{test
mapping} \emph{$\Lambda:F\rightarrow T$}, where the \emph{test space}
$T$ is assumed to be complete in order to ensure the boundedness
of $\Lambda^{-1}$. We denote
$\tilde{T}:=\Lambda\tilde{F}$. Test data from $T$ is discretized into finite-dimensional
test subspaces $T_{s}$ with a \emph{test discretization} mapping
\begin{equation}
\pi_{s}:T\rightarrow T_{s}\label{eq:test_discretization_definition}\end{equation}
the operator norm of which must be bounded by a constant, which is
independent of $s$. It follows that the operator norm of \[
\pi_{s}\Lambda L:U\rightarrow T_{s}\]
 is bounded similarly since \[
\left\Vert \pi_{s}\Lambda L\right\Vert _{U\rightarrow T_{s}}\leq\left\Vert \pi_{s}\right\Vert _{T\rightarrow T_{s}}\left\Vert \Lambda L\right\Vert _{U\rightarrow T}\]

In order to apply the error bound of Schaback's framework a number
of inequalities must be supplied. The first of these is the \emph{trial
space approximation property}\begin{equation}
\left\Vert u-\Pi_{r}u\right\Vert _{U}\leq\epsilon\left(r\right)\left\Vert u\right\Vert _{\tilde{U}}\mbox{ for all }u\in\tilde{U}.\label{eq:trial_space_approximation_property}\end{equation}
The second inequality is the test discretization's \emph{stability
condition} \begin{equation}
\left\Vert \Lambda Lu_{r}\right\Vert _{T}\leq\beta\left(s\right)\left\Vert \pi_{s}\Lambda Lu_{r}\right\Vert _{T_{s}}\mbox{ for all }u_{r}\in U_{r}.\label{eq:stability_condition}\end{equation}
If the \emph{stability factor} $\beta\left(s\right)$ grows as the
test discretization is refined, i.e., as $s$ decreases towards zero,
then the order of convergence in the final error bound~\eqref{eq:error_bound}
will be less than that provided by the trial space approximation property
\eqref{eq:trial_space_approximation_property}. When the stability
factor does not grow, the test discretization is called \emph{uniformly
stable}. The final inequality required by Schaback's framework involves
a numerical method capable of providing an approximate solution $u_{r,s}^{*}\in U_{r}$
which satisfies the \emph{numerical method approximation property}\begin{equation}
\left\Vert \pi_{s}\Lambda\left(Lu_{r,s}^{*}-f\right)\right\Vert _{T_{s}}\leq C\left\Vert \pi_{s}\Lambda L\right\Vert _{U\rightarrow T_{s}}\epsilon\left(r\right)\left\Vert u^{*}\right\Vert _{\tilde{U}}.\label{eq:numerical_method_approximation_property}\end{equation}
In particular, if the numerical method computes $u_{r,s}^{*}\in U_{r}$
which minimizes the left hand side of~\eqref{eq:numerical_method_approximation_property}
then the constant is at most one, since\begin{eqnarray}
\left\Vert \pi_{s}\Lambda\left(Lu_{r,s}^{*}-f\right)\right\Vert _{T_{s}} & \leq & \left\Vert \pi_{s}\Lambda L\left(\Pi_{r}u^{*}-u^{*}\right)\right\Vert _{T_{s}}\nonumber \\
 & \leq & \left\Vert \pi_{s}\Lambda L\right\Vert _{U\rightarrow T_{s}}\left\Vert \Pi_{r}u^{*}-u^{*}\right\Vert _{U}\nonumber \\
 & \leq & \left\Vert \pi_{s}\Lambda L\right\Vert _{U\rightarrow T_{s}}\epsilon\left(r\right)\left\Vert u^{*}\right\Vert _{\tilde{U}}.\label{eq:justification_for_nma_property}\end{eqnarray}
 
\begin{theorem}
\cite[Theorem 1]{schaback:2007:framework}\label{thm:error_bound}
Given the setting stated above, if the inequalities~\eqref{eq:trial_space_approximation_property},
\eqref{eq:stability_condition}, and~\eqref{eq:numerical_method_approximation_property}
are satisfied then the following error bound holds:
\begin{eqnarray*}
\left\Vert u^{*}-u_{r,s}^{*}\right\Vert _{U} & \leq & \left(1+\beta\left(s\right)\left\Vert \left(\Lambda L\right)^{-1}\right\Vert _{T\rightarrow U}\left\Vert \pi_{s}\Lambda L\right\Vert _{U\rightarrow T_{s}}\left(1+C\right)\right)\epsilon\left(r\right)\left\Vert u^{*}\right\Vert _{\tilde{U}}.\label{eq:error_bound}
\end{eqnarray*}

\end{theorem}
\begin{proof}
We have that\begin{eqnarray*}
\left\Vert u^{*}-u_{r,s}^{*}\right\Vert _{U} & \leq & \left\Vert u^{*}-\Pi_{r}u^{*}\right\Vert _{U}+\left\Vert \Pi_{r}u^{*}-u_{r,s}^{*}\right\Vert _{U}\\
 & \leq & \epsilon\left(r\right)\left\Vert u^{*}\right\Vert _{\tilde{U}}+\left\Vert \Pi_{r}u^{*}-u_{r,s}^{*}\right\Vert _{U}\end{eqnarray*}

\begin{eqnarray*}
\left\Vert \Pi_{r}u^{*}-u_{r,s}^{*}\right\Vert _{U} & \leq & \left\Vert \left(\Lambda L\right)^{-1}\right\Vert _{T\rightarrow U}\left\Vert \Lambda L\left(\Pi_{r}u^{*}-u_{r,s}^{*}\right)\right\Vert _{T}\\
 & \leq & \beta\left(s\right)\left\Vert \left(\Lambda L\right)^{-1}\right\Vert _{T\rightarrow U}\left\Vert \pi_{s}\Lambda L\left(\Pi_{r}u^{*}-u_{r,s}^{*}\right)\right\Vert _{T_{s}}\\
 & \leq & \beta\left(s\right)\left\Vert \left(\Lambda L\right)^{-1}\right\Vert _{T\rightarrow U}\left(\left\Vert \pi_{s}\Lambda L\left(\Pi_{r}u^{*}-u^{*}\right)\right\Vert _{T_{s}}\right.\\
 &      &  ~~~~~~~~~~~~~~~~~~~~~~~~+\left.\left\Vert \pi_{s}\Lambda L\left(u^{*}-u_{r,s}^{*}\right)\right\Vert _{T_{s}}\right)\\
 & = & \beta\left(s\right)\left\Vert \left(\Lambda L\right)^{-1}\right\Vert _{T\rightarrow U}\left\Vert \pi_{s}\Lambda L\right\Vert _{U\rightarrow T_{s}}\epsilon\left(r\right)\left\Vert u^{*}\right\Vert _{\tilde{U}}\left(1+C\right)\end{eqnarray*}

\begin{eqnarray*}
\left\Vert u^{*}-u_{r,s}^{*}\right\Vert _{U} & \leq & \epsilon\left(r\right)\left\Vert u^{*}\right\Vert _{\tilde{U}} \\
 &      & +\beta\left(s\right)\left\Vert \left(\Lambda L\right)^{-1}\right\Vert _{T\rightarrow U}\left\Vert \pi_{s}\Lambda L\right\Vert _{U\rightarrow T_{s}}\epsilon\left(r\right)\left\Vert u^{*}\right\Vert _{\tilde{U}}\left(1+C\right).
\end{eqnarray*}

\end{proof}
The stability condition~\eqref{eq:stability_condition} can be established
using an \emph{inverse estimate}\begin{equation}
\left\Vert u_{r}\right\Vert _{\tilde{U}}\leq\gamma\left(r\right)\left\Vert u_{r}\right\Vert _{U}\mbox{ for all }u_{r}\in U_{r},\label{eq:inverse_inequality}\end{equation}
 a \emph{sampling inequality} \begin{equation}
\left\Vert f\right\Vert _{T}\leq C\left(\alpha\left(s\right)\left\Vert f\right\Vert _{\tilde{T}}+\beta\left(s\right)\left\Vert \pi_{s}f\right\Vert _{T_{s}}\right)\mbox{ for all }f\in\tilde{T},\label{eq:sampling_inequality}\end{equation}
and ensuring that a \emph{fine enough test discretization} is chosen
such that \begin{equation}
C\alpha\left(s\right)\gamma\left(r\right)\left\Vert \Lambda L\right\Vert _{\tilde{U}\rightarrow\tilde{T}}\left\Vert \left(\Lambda L\right)^{-1}\right\Vert _{T\rightarrow U}\leq\frac{1}{2},\label{eq:fine_enough_test_discretization}\end{equation}
where $C$ is the constant appearing in~\eqref{eq:sampling_inequality}.
Typically, $\gamma\left(r\right)\rightarrow\infty$ as $r\rightarrow0$,
while $\alpha\left(s\right)\rightarrow0$ as $s\rightarrow0$.
\begin{proposition}
\label{thm:testing_stability_via_sampling_inequality}\cite[Theorem 2]{schaback:2007:framework}
If~\eqref{eq:inverse_inequality},~\eqref{eq:sampling_inequality},
and~\eqref{eq:fine_enough_test_discretization} hold then so does
\eqref{eq:stability_condition}.\end{proposition}
\begin{proof}
We have that\begin{eqnarray*}
\left\Vert \Lambda Lu_{r}\right\Vert _{T} & \leq & C\left(\alpha\left(s\right)\left\Vert \Lambda Lu_{r}\right\Vert _{\tilde{T}}+\beta\left(s\right)\left\Vert \pi_{s}\Lambda Lu_{r}\right\Vert _{T_{s}}\right)\\
 & \leq & C\left(\alpha\left(s\right)\left\Vert \Lambda L\right\Vert _{\tilde{U}\rightarrow\tilde{T}}\left\Vert u_{r}\right\Vert _{\tilde{U}}+\beta\left(s\right)\left\Vert \pi_{s}\Lambda Lu_{r}\right\Vert _{T_{s}}\right)\\
 & \leq & C\left(\alpha\left(s\right)\left\Vert \Lambda L\right\Vert _{\tilde{U}\rightarrow\tilde{T}}\gamma\left(r\right)\left\Vert u_{r}\right\Vert _{U}+\beta\left(s\right)\left\Vert \pi_{s}\Lambda Lu_{r}\right\Vert _{T_{s}}\right)\\
 & \leq & C\left(\alpha\left(s\right)\left\Vert \Lambda L\right\Vert _{\tilde{U}\rightarrow\tilde{T}}\gamma\left(r\right)\left\Vert \left(\Lambda L\right)^{-1}\right\Vert _{T\rightarrow U}\left\Vert \Lambda Lu_{r}\right\Vert _{T}\right.\\
 &      & ~~~~\left.+\beta\left(s\right)\left\Vert \pi_{s}\Lambda Lu_{r}\right\Vert _{T_{s}}\right)\\
 & \leq & \frac{1}{2}\left\Vert \Lambda Lu_{r}\right\Vert _{T}+C\beta\left(s\right)\left\Vert \pi_{s}\Lambda Lu_{r}\right\Vert _{T_{s}}\end{eqnarray*}
and the result follows by incorporating the constant $2C$ in $\beta\left(s\right)$.  
\end{proof}

\subsection{\label{sec:PoissonProblem}Convergence Results for the Poisson Problem}

We consider the example from~\cite[Section 4.1]{schaback:2007:framework},
a Poisson problem with mixed, inhomogeneous boundary data:
let $\Omega$ be a bounded domain in $\mathbb{R}^{d}$ with a Lipschitz-continuous
boundary.  We denote $\Omega_{1}:=\Omega$, $\Omega_{2}:=\Gamma^{D}\subset\partial\Omega$,
and $\Omega_{3}=\Gamma^{N}\subset\partial\Omega$ so that the dimension of each domain
is given by $n_{1}=n$, and $n_{2},n_{3}=n-1$. 
Let $m,\tilde{m}$ be nonnegative real numbers such that $\tilde{m}-m\in\mathbb{N}$, and
\begin{eqnarray}
\left(m_1, m_2, m_3\right) & := & \left(m,m+3/2,m+1/2\right) \nonumber \\
U & := & H^{m+2}\left(\Omega\right)\nonumber \\
F & := & F^1 \times F^2 \times F^3 \nonumber \\
  & :=  & H^{m_1}\left(\Omega_1\right)\times H^{m_2}\left(\Omega_2\right)\times H^{m_3}\left(\Omega_3\right)\nonumber \\
Lu & := & \left(-\Delta u,u|_{\Gamma^{D}},\frac{\partial u}{\partial n}|_{\Gamma^{N}}\right), \label{eq:Poisson_problem}\end{eqnarray}
with analogous definitions made for $\left(\tilde{m}_1,\tilde{m}_2,\tilde{m}_3\right)$, $\tilde{U}$, and $\tilde{F}$.
With the space $F$ equipped with the norm
$\left\Vert \cdot\right\Vert _{F}^{2}:=\left\Vert \cdot\right\Vert _{F^{1}}^{2}+\left\Vert \cdot\right\Vert _{F^{2}}^{2}+\left\Vert \cdot\right\Vert _{F^{3}}^{2}$,
it follows that the linear operator $L$, as defined above, is continuously invertible either as $L:U\rightarrow F$ or $L:\tilde{U}\rightarrow \tilde{F}$.

We assume that the solution comes from $\tilde{U}$ and that the trial space $U_{r}$
is chosen such that the trial space approximation property~\eqref{eq:trial_space_approximation_property}
holds with $\epsilon\left(r\right)=O\left(r^{\tilde{m}-m}\right)$,
a property satisfied by kernel-based meshless trial spaces, c.f. Narcowich et al.~\cite{narc_ward_wend:scattered_zeros:2005,narco:2006:error_bound}, and 
finite-element trial spaces~\cite[Theorem 4.5.11]{brenner_scott}.
We also assume that the inverse estimate~\eqref{eq:inverse_inequality}
holds with $\gamma\left(r\right)=O\left(r^{m-\tilde{m}}\right)$,
as is the case for finite-element trial spaces~\cite[Theorem 4.4.20]{brenner_scott}.
Obtaining an inverse estimate with the expected factor $\gamma\left(r\right)=O\left(r^{m-\tilde{m}}\right)$
appears to be an open problem for kernel-based meshless trial spaces.
Narcowich et al.~\cite{narco:2006:error_bound} provide an inverse
estimate with the expected factor for the case of Sobolev spaces over
$\mathbb{R}^{n}$. Both Schaback and Wendland~\cite{schaback-wendland:2002-1},
and Duan~\cite{duan:inverse-estimate:2008} provide inverse estimates
for Sobolev spaces over a domain. Unfortunately, the factor involved
in these inverse estimates are worse than the finite-element case.
Further progress on this problem is expected to be reported in the
thesis of Rieger~\cite{rieger:2008:thesis}.

We consider the case of strong testing here, which means that the
test mapping $\Lambda:F\rightarrow T$ is just the identity mapping and
that each test space $T^k$ coincides with the corresponding data space $F^k$.
Weak testing is also possible, in which case the test functionals integrate functions in $F^k$
against test functions, resulting in the test data in each $T^k$ acquiring additional smoothness.
This is discussed in detail by Schaback~\cite{schaback:2007:framework,schaback:2007:weak_testing}.
Each domain is discretized onto finite
subsets $Y_{s}^{k}\subset\Omega_{k}$, with the same fill distance
$s=\delta\left(Y_{s}^{k},\Omega_{k}\right)$. Furthermore, they are
assumed to satisfy the property that $\#Y_{s}^{k}$ is bounded by
$s^{-n_{k}}$ up to a constant, as is the case for domain discretization
with a uniformly bounded mesh ratio~\cite{schaback:2007:framework}.
We define discrete test spaces \[
T_{s}^{k}:=\mathbb{R}^{\#\left\{ \alpha:\left|\alpha\right|\leq\mu_{k}\right\} \cdot\#Y_{s}^{k}}\]
equipped with a norm \begin{equation}
\left\Vert \cdot\right\Vert _{T_{s}^{k}}:=s^{n_{k}/2}\left\Vert \cdot\right\Vert _{2}\label{eq:discrete_test_norm}\end{equation}
and a test discretization $\pi_{s}^{k}:T^{k}\rightarrow T_{s}^{k}$
\begin{equation}
\pi_{s}^{k}f_{k}:=\prod_{\left|\alpha\right|\leq\mu_{k}}\left.\partial^{\alpha}f_{k}\right|_{Y_{s}^{k}}\mbox{ for all }f_{k}\in T^{k} \label{eq:test_discretization_component_def}\end{equation}
where $\mu_{k}$ is an integer such that 
$m_{k}-\mu_{k}-n_k/2 > 0$, and furthermore this difference is independent of $k$. 
The discrete test space $T_s:=T_s^1\times T_s^2\times T_s^3$ is defined and equipped with a norm, analogously to $F$ and $T$. The test space $T$
is then equipped with a test discretization $\pi_{s}:T\rightarrow T_{s}$
defined by \begin{equation}
\pi_{s}f:=\left(\pi_{s}^{1}f_{1},\pi_{s}^{2}f_{2},\pi_{s}^{3}f_{3}\right)\mbox{ for all }f=\left(f_{1},f_{2},f_{3}\right)\in T,
\label{eq:test_discretization_def}
\end{equation}
\begin{proposition}
If for each $k$, $m_k-\mu_k-n_k/2>0$  then $\pi_{s}:T\rightarrow T_{s}$ is well-defined
and the operator norm $\left\Vert \pi_{s}\right\Vert _{T\rightarrow T_{s}}$
is bounded independently of $s$.\end{proposition}
\begin{proof}
Suppose $f=\left(f_{1},f_{2},f_{3}\right)\in T$.
Since $m_{k}-\mu_{k}>n_k/2$ we have from the Sobolev embedding theorem
that $T^{k}\hookrightarrow C^{\mu_{k}}\left(\overline{\Omega_{k}}\right)$
and therefore the test discretization is both well-defined and there
exists some constant independent of $f=\left(f_{1},f_{2},f_{3}\right)$
such that for each $f_{k}$,\[
\left\Vert f_{k}\right\Vert _{C^{\mu_{k}}\left(\overline{\Omega_{k}}\right)}\leq C\left\Vert f_{k}\right\Vert _{T^{k}}.\]
Since $\#Y_{s}^{k}$ is bounded by $s^{-n_{k}}$ up to some constant
which is independent of $s$, it follows that \begin{eqnarray*}
\left\Vert \pi_{s}f\right\Vert _{T_{s}}^{2} & = & \sum_{k=1}^{3}\left\Vert \pi_{s}^{k}f_{k}\right\Vert _{T_{s}^{k}}^{2}\\
 & = & \sum_{k=1}^{3}s^{n_{k}}\sum_{x\in Y_{s}^{k}}\sum_{\left|\alpha\right|\leq\mu_{k}}\left|\partial^{\alpha}f\left(x\right)\right|^{2}\\
 & \leq & \sum_{k=1}^{3}s^{n_{k}}\left\Vert f\right\Vert _{C^{\mu_{k}}\left(\overline{\Omega_{k}}\right)}^{2}\#\left\{ \alpha:\left|\alpha\right|\leq\mu_{k}\right\} \cdot\#Y_{s}^{k}\\
 & \leq & C\sum_{k=1}^{3}\left\Vert f\right\Vert _{T^{k}}^{2}=C\left\Vert f\right\Vert _{T}^{2}.\end{eqnarray*}
\end{proof}

\begin{proposition}
\label{pro:full_sampling_inequality}There exists a constant $s_{0}$
such that for all $s\leq s_{0}$ a  sampling inequality~\eqref{eq:sampling_inequality}
holds with a constant $C$ for the test space $T$ and test discretization $\pi_{s}:T\rightarrow T_{s}$
with $\alpha\left(s\right):=s^{\tilde{m}-m}$, and $\beta\left(s\right):=s^{\mu_1-m_1}=s^{\mu_2-1/2-m_2}=s^{\mu_3-1/2-m_3}$.\end{proposition}
\begin{proof}
From Corollary~\ref{cor:specialized_sampling_inequality} and~\eqref{eq:discrete_test_norm}, it follows that for each $k$
there exist constants $C_{k}$ and $s_{k}$ such that for $s\leq s_0 := \min\left(1,s_1,s_2,s_3\right)$,
\begin{eqnarray*}
\left\Vert f_{k}\right\Vert _{T^{k}} & \leq & C_{k}\left(\alpha\left(s\right)\left\Vert f_{k}\right\Vert _{\tilde{T}^{k}}+s^{\mu_{k}-m_{k}}\left\Vert \pi_{s}^{k}f_{k}\right\Vert _{T^{k}_{s}}\right)\\
 & \leq & C_{k}\left(\alpha\left(s\right)\left\Vert f_{k}\right\Vert _{\tilde{T}^{k}}+s^{\mu_{1}-m_{1}}\left\Vert \pi_{s}^{k}f_{k}\right\Vert _{T^{k}_{s}}\right)\end{eqnarray*}
since $\mu_{1}-m_{1}\leq\mu_{2}-m_{2}=\mu_{3}-m_{3}$.  The result then follows with a constant $C$ by combining the preceding inequalities.
\end{proof}
We assume that the $s$ is sufficiently
small to satisfy the requirements of Proposition~\ref{pro:full_sampling_inequality}
and~\eqref{eq:fine_enough_test_discretization}. Even in the fractional
case, the sampling inequality introduced here provides $\alpha\left(s\right)=s^{\tilde{m}-m}$
which shrinks as rapidly as the expected inverse estimate factor $\gamma\left(r\right)=r^{m-\tilde{m}}$
grows and thus $s$ and $r$ can be kept proportional. This is in
contrast to previous sampling inequalities which necessarily introduce a factor $\alpha\left(s\right)=s^{\tilde{m}-\left\lceil m\right\rceil }$
when bounding fractional order Sobolev norms, requiring the test discretization
to be refined more rapidly than the trial discretization and thus
diminishing the order of convergence by $\left\lceil m\right\rceil -m$.
If the function $u_{r,s}^{*}\in U_{r}$ which minimizes the left hand
side of~\eqref{eq:numerical_method_approximation_property} has been
computed, then Schaback's framework provides the error bound~\eqref{eq:error_bound}
with constant $C=1$. The order of convergence established by this error bound, in terms of both the
trial and test discretization, is then given by $\beta\left(h\right)\epsilon\left(h\right)$ and
satisfies \[
\beta\left(h\right)\epsilon\left(h\right)=O\left(h^{\left(\tilde{m}-m\right)+\left(\mu_1-m_1\right)}\right).\]

Table~\ref{fig:TwoThreeDConvergenceRates}
states particular convergence results using various Sobolev norms and
test discretizations in two- and three-dimensions. These results show that, for convergence in higher order norms,
the highest order of convergence is obtained using a higher order test discretization introduced here.

We note that to obtain a uniformly stable test discretization with
$\beta\left(s\right)=1$ would require choosing the order $\mu_{k}$
of each test discretization to be equal to that of the order $m_{k}$
of the test space. Unfortunately, this does not seem to be possible:
in order for the test discretizations operator norm to be bounded
independently of $s$, the order of each test space $T^{k}$ is required
to be greater than that of the test discretization $\mu_{k}$ by at
least $n_{k}/2$. It follows that the order of convergence, in terms
of both the trial and test discretization,  provided by this modified
formulation of Schaback's framework is always less than that of the
trial space approximation property. Another consequence is that the
order of $U$ must be at least $2+n/2$, and therefore convergence
in the $L^{2}$ norm can only be concluded suboptimally from convergence results in higher order Sobolev norms, using strong testing in this modified formulation
of Schaback's framework.

\begin{table}
\begin{tabular}{|c|c|c|c|c|}
\hline 
$U=H^{m+2}\left(\Omega\right)=$ & $H^{0}\left(\Omega\right)$ & $H^{4}\left(\Omega\right)$ & $H^{5}\left(\Omega\right)$ & $H^{6}\left(\Omega\right)$\tabularnewline
\hline
\hline 
$\mu_{1} = 0$ & None & $\tilde{m}-m-2$ & $\tilde{m}-m-3$ & $\tilde{m}-m-4$\tabularnewline
\hline 
$\mu_{1} = 1$ & None & None & $\tilde{m}-m-2$ & $\tilde{m}-m-3$\tabularnewline
\hline 
$\mu_{1} = 2$ & None & None & None & $\tilde{m}-m-2$\tabularnewline
\hline 
$\mu_{1} = 3$ & None & None & None & None\tabularnewline
\hline
\end{tabular}

\caption{\label{fig:TwoThreeDConvergenceRates}Order of convergence in various
Sobolev norms established by a modified formulation of Schaback's
framework, using trial spaces with optimal properties and strong testing
with various order test discretizations to solve two- or three-dimensional
Poisson problems.}

\end{table}

\section{Conclusions}

We have further generalized the sampling inequalities of Arcang\'{e}li
et al.~\cite{arcangeli:2007:sampling_inequality}, Madych~\cite{madych:2006:sampling_inequality},
and Wendland and Rieger~\cite{wendland:2005:sampling_inequality},
to optimally bound fractional order Sobolev semi-norms, and to incorporate
higher order data into the bound. When used in a modified formulation
of Schaback's framework to prove convergence rates for unsymmetric
meshless methods this new sampling inequality has two benefits:
\begin{enumerate}
\item It results in more optimal estimates for problems involving fractional
order Sobolev spaces, particularly by providing a more optimal constant
$\alpha\left(s\right)$.
\item For convergence in higher order Sobolev norms, higher order results are obtained using a higher order test discretization in comparison to the zero order test discretization.
\end{enumerate}
The zero order test discretization has been widely employed
in practice, and corresponds to what is usually called Kansa's method or unsymmetric collocation.
On the other hand higher order testing has not, and its value in practical
applications requiring convergence in stronger norms is an open question
worthy of further study.

\bibliographystyle{amsplain}
\bibliography{general_sampling_inequality}

\providecommand{\bysame}{\leavevmode\hbox to3em{\hrulefill}\thinspace}
\providecommand{\MR}{\relax\ifhmode\unskip\space\fi MR }
\providecommand{\MRhref}[2]{%
  \href{http://www.ams.org/mathscinet-getitem?mr=#1}{#2}
}
\providecommand{\href}[2]{#2}
\begin{thebibliography}{10}

\bibitem{arcangeli:2007:sampling_inequality}
R\'emi Arcang\'eli, Mar\'ia Cruz~L\'opez de~Silanes, and Juan~Jos\'e Torrens,
  \emph{An extension of a bound for functions in {S}obolev spaces, with
  applications to (m,s)-spline interpolation and smoothing}, Numerische
  Mathematik \textbf{107} (2007), 181--211.

\bibitem{book_with_fractional_paper:2001}
Jean Bourgain, Haim Brezis, and Petru Mironescu, \emph{Optimal control and
  partial differential equations, in honour of {P}rofessor {A}lain
  {B}ensoussan's 60th birthday}, ch.~Another Look at Sobolev Spaces,
  pp.~439--455, IOS Press, 2001.

\bibitem{brenner_scott}
Susanne~C. Brenner and L.~Ridgway Scott, \emph{The mathematical theory of
  finite element methods}, third ed., Texts in Applied Mathematics, vol.~15,
  Springer, New York, 2008.

\bibitem{duan:inverse-estimate:2008}
Yong Duan, \emph{A note on the meshless method using radial basis functions},
  Computers \& Mathematics with Applications \textbf{55} (2008), 66--75.

\bibitem{evans:pde}
Lawrence Evans, \emph{Partial differential equations}, AMS, 2002.

\bibitem{madych:2006:sampling_inequality}
W.R. Madych, \emph{An estimate for multivariate interpolation ii}, Journal of
  approximation theory \textbf{142} (2006), 116--128.

\bibitem{narc_ward_wend:scattered_zeros:2005}
Francis~J. Narcowich, Joseph~D. Ward, and Holger Wendland, \emph{Sobolev bounds
  on functions with scattered zeros, with applications to radial basis function
  surface fitting}, Math. Comp. \textbf{74} (2005), no.~250, 743--763
  (electronic).

\bibitem{narco:2006:error_bound}
Francis~J. Narcowich, Joseph~D. Ward, and Holger Wendland, \emph{Sobolev error
  estimates and a {B}ernstein inequality for scattered data interpolation via
  radial basis functions}, Constructive Approximation \textbf{24} (2006),
  175--186.

\bibitem{rieger:2008:thesis}
Christian Rieger, \emph{Sampling inequalities and applications}, Ph.D. thesis,
  G\"ottingen, 2008.

\bibitem{rieger:2006:exponential_sampling_inequality}
Christian Rieger and Barbara Zwicknagl, \emph{Sampling inequalities for
  infinitely smooth functions, with applications to interpolation and machine
  learning}, 2008.

\bibitem{schaback-wendland:2002-1}
R.~Schaback and H.~Wendland, \emph{Inverse and saturation theorems for radial
  basis function interpolation}, Math. Comp. \textbf{71} (2002), no.~238,
  669--681 (electronic). \MR{MR1885620 (2003a:41018)}

\bibitem{schaback:2007:framework}
Robert Schaback, \emph{Unsymmetric meshless methods for operator equations},
  Preprint G\"ottingen, 2006.

\bibitem{schaback:2007:meshless_collocation_convergence}
\bysame, \emph{Convergence of unsymmetric kernel-based meshless collocation
  methods}, SIAM Journal on Numerical Analysis \textbf{45} (2007), no.~1,
  333--351.

\bibitem{schaback:2007:weak_testing}
\bysame, \emph{Recovery of functions from weak data using unsymmetric meshless
  kernel-based methods}, To appear in Applied Numerical Mathematics, 2007.

\bibitem{wendland:2005:sampling_inequality}
H.~Wendland and C.~Rieger, \emph{Approximate interpolation with applications to
  selecting smoothing parameters}, Numerische Mathematik \textbf{101} (2005),
  729--748.

\end{thebibliography}

\end{document}